\documentclass{amsart}

\usepackage{graphicx}

\newtheorem{theorem}{Theorem}[section]
\newtheorem*{theorem A}{Theorem A}
\newtheorem*{theorem B}{N\"olker's Theorem}
\newtheorem{lemma}{Lemma}[section]

\newtheorem{notation}{Notation}[section]

\newtheorem*{question}{Question}

\newtheorem {conjecture}{Conjecture}
\theoremstyle{remark}
\newtheorem{remark}{Remark}[section]
\theoremstyle{remark}

\theoremstyle{definition}

\newtheorem{definition}{Definition}[section]

\numberwithin{equation}{section}
\def\({\left ( }
\def\){\right )}
\def\<{\left < }
\def\>{\right >}

\begin{document}

%\noindent {\sc International Electronic Journal of Geometry}

%\noindent {\sc \small Volume 1  No. ? pp. 000--000 (2008) \copyright
%IEJG}

\vspace{2cm}

\title {A Lower Bound for the Mahler Volume of Symmetric Convex Sets}

%    Information for first author
\author{Yashar Memarian}   
\address{Last Institution: Department of Mathematics\\
    University of Notre-Dame}
\email{y.memarian@gmail.com}

\thanks{I like to thank F.Barthe, M.Fradelizi and R.Vershynin for their useful comments about this paper.}


%\date{October 15, 2015}


\begin{abstract}
The goal of this paper is to present a lower bound for the Mahler volume of at least 4-dimensional symmetric convex bodies. We define a computable dimension dependent constant through a 2-dimensional variational (max-min) procedure and demonstrate that the Mahler volume of every (at least 4-dimensional) symmetric convex body is greater than a (simple) function of this constant. Similar to the proof of Gromov's Waist of the Sphere Theorem in \cite{grwst}, our result is proved via localisation-type arguments obtained from a suitable measurable partition (or partitions) of the canonical sphere.
\end{abstract}
\maketitle

\section{Introduction}
The Mahler Conjecture is one of those amusing conjectures that have been floating around mathematical literature for decades. It originates (as the name indicates) from Mahler's studies at the end of the Thirties. It is one of those mathematical conjectures that everyone with a minimal knowledge in  mathematics can understand. The idea is simple: consider (symmetric) convex sets in a fixed dimension and ask what the \emph{most} and \emph{least} (symmetric) ones are. To be able to give a mathematical flavor to this question, Mahler introduced a \emph{functional} on the class of (symmetric) convex sets, and  was interested in minimizing and maximizing this functional. The functional is now known as the \emph{Mahler volume} of a convex set. To understand this functional one should first know the definition of the polar of a (symmetric) convex body :
\begin{definition}[Polar of a Convex Body]
Let $K$ be a convex body in $\mathbb{R}^n$. The polar of $K$ denoted by $K^{\circ}$ is the following:
\begin{eqnarray*}
K^{\circ}=\{x\in\mathbb{R}^n\vert \vert x.y\vert\,\leq 1, \forall y\in K\},
\end{eqnarray*}
where $.$ stands for the inner product associated to the canonical Euclidean structure of $\mathbb{R}^n$.
\end{definition}
It is not hard to verify that indeed $K^{\circ}$ is a convex body (should $K$ be convex itself). Additionally, if we suppose $K$ to be symmetric with respect to the origin of $\mathbb{R}^n$, $K^{\circ}$ satisfies the same property. 

The Mahler volume of a (symmetric) convex body is defined by :
\begin{eqnarray*}
M(K)=vol_n(K)vol_n(K^{\circ}),
\end{eqnarray*}
and the interesting question is now to find:
\begin{eqnarray*}
\min_{K}M(K),
\end{eqnarray*}
where $K$ runs over the class of (symmetric) convex bodies.

Mahler raised this question in $1938-1939$ in \cite{mah1} and \cite{mah2}. He succeeded to completely answer his own question for the $2$-dimensional bodies. He realised that if one minimises the Mahler volume over the class of convex bodies, the minimum of the Mahler volume is attained for the simplex (i.e. the convex hull of three non-aligned points). When one seeks the minimum over symmetric convex bodies, the minimum is attained on the $2$-dimensional cube (see \cite{meyer} for another proof of the $2$-dimensional Mahler question). After this discovery, he naturally generalized the question for higher dimensions. Here, we are only interested in symmetric convex bodies:
\begin{conjecture}[Mahler Conjecture] \label{mal}
Let $K$ be a symmetric convex body in $\mathbb{R}^n$, then
\begin{eqnarray*}
vol_n(K)vol_n(K^{\circ}) &\geq& vol_n(I_n)vol_n(I_n^{\circ})\\
                          &=&\frac{4^n}{\Gamma(n+1)}.
\end{eqnarray*}
where $\Gamma(n)=(n-1)!$ and is the (well-known) Gamma function :
\begin{eqnarray*}
\Gamma(t)=\int_{0}^{\infty}x^{t-1}e^{-x}dx,
\end{eqnarray*}
and where $I_n$ is the $n$-dimensional unit cube.
\end{conjecture}

Mahler did not just conjecture the above. He was also interested in an upper bound for the Mahler volume of (symmetric) convex bodies. He conjectured that the upper bound is achieved by the unit ball. The upper bound case for the Mahler question has been completely answered and is known as the Blaschke-Santalo inequality:
\begin{theorem}[Blaschke-Santalo Inequality] \label{bl}
For every (symmetric) convex body $K$ in $\mathbb{R}^n$ where $n\geq 2$, we have:
\begin{eqnarray*}
vol_n(K)vol_n(K^{\circ})&\leq& vol_n(B_n(0,1))vol_n(B_n(0,1)^{\circ})\\
                    &=& vol_n(B_n(0,1))^2 \\
                    &=& \frac{\pi^{n/2}}{\Gamma(\frac{n}{2}+1)},
\end{eqnarray*}
where $B_n(0,1)$ is the $n$-dimensional unit ball in $\mathbb{R}^n$.
\end{theorem}
Theorem \ref{bl} (which was one part of Mahler's question) has been proved by Blaschke for the $3$-dimensional case in \cite{bl} and \cite{bl2} and in full generality by Santalo in \cite{santal}. There are also other proofs to this theorem : in \cite{meypaj} the authors prove this theorem using Steiner symmetrisation. In \cite{sr} the author proves that the equality case in Theorem \ref{bl} is obtained only for ellipsoids. In \cite{leh1} the author proves a functional inequality version of this inequality using a partition-type method. One can also consult Tao's blog post \cite{tao1} and \cite{lut} for more information on this inequality.

The lower bound for the Mahler volume of (symmetric) convex bodies however is still open for $n\geq 3$. Naturally, since Mahler raised Conjecture \ref{mal} many developments have been made concerning this conjecture. There are many partial results regarding this question which can be settled in different ways.
One difficulty with Conjecture \ref{mal} is that the minimizer is not \emph{unique}. Indeed the Mahler volume is invariant under affine transformations of $\mathbb{R}^n$. This complicates finding this minimizer. There is a very nice note about the complications within Conjecture \ref{mal} in \cite{tao2}. Although one can (easily) deduce that a minimizer for the Mahler volume exists (Mahler himself was able to demonstrate this fact), not much was known about the \emph{shape} of such a minimizer until very recently. In \cite{stan}, the author demonstrates that the boundary of the minimizer can not be positively curved and of class $C^2$ everywhere. In \cite{reisut}, the authors prove a result suggesting that the minimizer of the Mahler volume is indeed a polytope.

Though the Mahler Conjecture suggests that the global minimizer for the Mahler volume is the cube, one could also seek the local minimizer of this functional. In \cite{naz}, the authors demonstrate that the (unit) cube is indeed a local minimizer for the Mahler functional.

As proving Conjecture \ref{mal} in full generality has revealed many complications, a possible way to prove partial results related to this question would involve trying to prove the conjecture for more restricted (symmetric) convex bodies. In \cite{sr} and \cite{reis3}, the authors prove the conjecture if one restricts themself to the class of unit balls of Banach spaces with a $1$-unconditional basis. In \cite{reis2} and \cite{reis}, the author proves the conjecture for the Zonoids. A simplified proof of this fact can be found in \cite{gord}. In \cite{barthe}, the authors prove the conjecture for symmetric convex bodies with \emph{many} symmetries. See also \cite{lop} for a partial result concerning Conjecture \ref{mal}.

Another direction one could take to study Conjecture \ref{mal} would be by trying to find non-necessarily sharp lower bounds for the Mahler volume of (symmetric) convex bodies. Usually this is done by finding a constant $c(n)$ (depending on the dimension) such that the Mahler volume of \emph{every} (symmetric) convex body is larger than $c(n)M(B_n(0,1))$. The most trivial constant one can easily obtain is by using the famous John's Lemma (see \cite{john}) to obtain $c(n)=n^{-n/2}$. The first breakthrough is due to Bourgain-Milman in \cite{bourg} where the authors prove the existence of a universal constant $c>0$ such that $c(n)=c^n$. In \cite{kuper1} the author discovers another value for $c(n)$ and (to my  knowledge), in \cite{kuper} the same author proves the best known constant for $c(n)$.

Although Conjecture \ref{mal} can be categorized as a geometric problem, it has been demonstrated that much information about it can be obtained by using functional inequalities and analysis. For example, the Legendre transformation has had great importance in studying this conjecture (see \cite{meyreis}, \cite{fradmey}, \cite{white} and \cite{groe}).

Numerous surveys and books on Convex Geometry (in which at least a few chapters are dedicated to the volume of polar bodies and/or the Mahler volume) can be consulted in: \cite{ball}, \cite{versh}, \cite{pis}, \cite{sch}, \cite{heuze},\cite{hu} and \cite{gard}.

The Mahler Conjecture has shown importance in other areas of mathematics as well. For example, in \cite{ball1} the author shows the connection between Conjecture \ref{mal} to wavelets, and recently in \cite{shiri} the authors connect Conjecture \ref{mal} to questions in symplectic geometry.

The main result of this paper presents a lower bound for the Mahler volume of symmetric convex bodies when the dimension of the sets satisfies $n\geq 4$. This lower bound is obtained from a variational procedure in $\mathbb{R}^2$. As we observed above, the Mahler Conjecture itself can be seen as the variational problem of minimizing the Mahler functional over the class of symmetric convex sets. The variational procedure with which we obtain this lower bound is rather non-trivial, but we believe it demonstrates the complications hidden in the Mahler Conjecture itself.

Before announcing the main theorem of this paper, we shall need a few definitions. Since the lower bound is presented through a procedure in $\mathbb{R}^2$, we should also give the necessary definitions in $\mathbb{R}^2$. The sets with which we work belong to a certain class of symmetric convex sets. For every $n\geq 4$ the class $\mathcal{S}(1,\sqrt{n+1})$ denotes the class of symmetric convex bodies $M$ in $\mathbb{R}^2$ such that
\begin{eqnarray*}
B_2(0,1)\subseteq M\subseteq B_2(0,\sqrt{n+1}).
\end{eqnarray*}
We no longer will be working with the Lebesgue measure- instead we shall define a class of measures which have an anisotropic density (i.e. the density of the measures will not be radial functions- they bear weight on the circle when writing the measure in polar coordinates).

\begin{definition}[The measures $\mu_{2,\theta}$] \label{mes}

Fix the Cartesian coordinates $(x,y)$ in $\mathbb{R}^2$. The measure $\mu_2$ is the measure defined as $r^n g(t)dr\wedge dt$ in polar coordinates of $\mathbb{R}^2$ where $g(t)=\vert\cos(t)^{n-1}\vert$. In Cartesian coordinates $(x,y)$, the measure $\mu_2$ is written:
\begin{eqnarray*}
\mu_2=\vert y\vert^{n-1}dx\,dy.
\end{eqnarray*}
Here, the end point of the unit vector of the $y$-axis coincides with the maximum point of the function $g$.
Given $\theta\in[0,\pi]$, we define the measure $\mu_{2,\theta}=r^n g(t+\theta)dr\wedge dt$. 
%The Cartesian coordinates associated to the measure $\mu_{2,\theta}$ coincides with the measure $\mu_2$.
\end{definition}
We will often work on conical subsets in $\mathbb{R}^2$:
\begin{notation}
We denote points with both non-negative coordinates of the unit circle by $\mathbb{S}^1_{+}$. If $I$ is a (geodesic) segment of $\mathbb{S}^1_{+}$, the cone $C(I)$ over $I$ will be the convex cone generated by $I$. The class of connected closed subsets of the unit circle (or the class of geodesic segments of the unit circle) will be denoted by $\mathcal{A}$.
\end{notation}

\begin{definition}[Constant $\alpha(n)$] \label{deff}
Let $n\geq 4$. Let $S$ be a $2$-dimensional convex set symmetric with respect to the origin of $\mathbb{R}^2$ and belonging to $\mathcal{S}(1,\sqrt{n+1})$. Define
\begin{eqnarray*}
\alpha(n,\theta,I,S)=\frac{\mu_{2,\theta}(C(I)\cap S)\mu_{2,\theta}(C(I)\cap S^{\circ})}{\big(\displaystyle\int_{I}g(\theta,t)dt\big)^2},
\end{eqnarray*}
where $\mu_{2,\theta}$ is the measure defined in \ref{mes} and $g(\theta,t)$ is the density of this measure in polar coordinates. Define
\begin{eqnarray*}
\alpha(n,S)=\max_{\theta\in[0,\pi]}\big(\min_{I\in \mathcal{A}}\alpha(n,\theta,I,S)\big).
\end{eqnarray*}
And at last, define:
\begin{eqnarray*}
\alpha(n)=\min_{S}\big(\alpha(n,S)\big),
\end{eqnarray*}
where $S$ runs over the class of $2$-dimensional symmetric convex sets in $\mathcal{S}(1,\sqrt{n+1})$. 
\end{definition}
For clarity, we shall explain how the constant $\alpha(n)$ is defined:
\begin{itemize}
\item First: fix a symmetric convex set $S$ in $\mathbb{R}^2$ and in $\mathcal{S}(1,\sqrt{n+1})$ where the Cartesian coordinates are also fixed.
\item Second step: fix a real $\theta\in[0,\pi]$.
\item With respect to the chosen $\theta$, consider the measure $\mu_{2,\theta}$ being defined on (the Borel subsets of) $\mathbb{R}^2$.
\item Find a geodesic segment $I\subset \mathbb{S}^1_{+}$ and consider the convex cone $C(I)$ generated by $I$.
\item We are now prepared to calculate the following:
\begin{eqnarray*}
\frac{\mu_{2,\theta}(C(I)\cap S)\mu_{2,\theta}(C(I)\cap S^{\circ})}{\big(\displaystyle\int_{I}g(\theta,t)dt\big)^2}.
\end{eqnarray*}
This will be the value $\alpha(n,\theta,I,S)$.
\item Now we are ready to follow a max-min procedure: by fixing $\theta$, take the min of the values $\alpha(n,\theta,I,S)$ where $I$ will pass over geodesic segments on the half-circle. Then we allow $\theta$ to move in $[0,\pi]$ (which will change the measures at every step) and calculate the maximum of $\alpha(n,\theta,I,S)$. This value will be the $\alpha(n,S)$.
\item All that remains now is to take $\alpha(n)$ as the minimum of the above over every symmetric convex set in $\mathcal{S}(1,\sqrt{n+1})$.
\end{itemize}

\begin{remark}
There is an \emph{easier} way to \emph{calculate} $\alpha(n,S)$ by introducing a variable $\phi$. Since we \emph{rotate} the measure $\mu_2$, instead of letting the intervals $I$ pass over all the connected intervals of the half-circle, we could consider the following:
\begin{eqnarray*}
\alpha_1(n,S)=\max_{\theta}\big(\min_{\phi}\alpha(n,\theta,[0,\phi],S)\big),
\end{eqnarray*}
where $\phi\in[0,\pi]$ and $[0,\phi]$ is an interval on the half-circle where the point $0$ coincides with the point $(1,0)$ in the Cartesian coordinates and the point $\phi$ coincides with the point $(\cos(\phi),\sin(\phi))$ on the half-circle. It is not hard to verify that in fact:
\begin{eqnarray*}
\alpha(n,S)=\alpha_1(n,S).
\end{eqnarray*}
\end{remark}
Before going any further, we illustrate that the constant $\alpha(n)$ has a non-zero lower bound, otherwise its definition would be meaningless:
\begin{lemma} \label{nontriv}
Let $\alpha(n)$ be defined as in \ref{deff}, then we have :
\begin{eqnarray*}
\alpha(n)\geq \frac{1}{(n+1)^2 (n+1)^{\frac{n+1}{2}}}.
\end{eqnarray*}
\end{lemma}

\begin{proof}

By definition, every $M\in\mathcal{S}(1,\sqrt{n+1})$ satisfies :
\begin{equation} \label{eqn:str}
B_2(0,1)\subseteq M\subseteq B_2(0,\sqrt{n+1}).
\end{equation}
If $K$ and $L$ are two convex sets such that $K\subset L$, then $L^{\circ}\subset K^{\circ}$. This shows that according to (\ref{eqn:str}) we have : 
\begin{eqnarray*}
B_2(0,\frac{1}{\sqrt{n+1}})\subseteq M^{\circ}.
\end{eqnarray*}
Let $I$ be any (geodesic) segment of $\mathbb{S}^1_{+}$ and let $C(I)$ be the convex cone generated by $I$. Let $\mu_{2,\theta}$ be a measure as in definition \ref{mes}. Therefore we have:
\begin{eqnarray*}
\alpha(n,\theta,I,K)&=&\frac{\mu_{2,\theta}(C(I)\cap K)\mu_{2,\theta}(C(I)\cap K^{\circ})}{\big(\displaystyle\int_{I}g(\theta,t)dt\big)^2}\\
                    &\geq& \frac{\mu_{2,\theta}(C(I)\cap B_2(0,1))\mu_{2,\theta}(C(I)\cap B_2(0,1/\sqrt{n+1}))}{\big(\displaystyle\int_{I}g(\theta,t)dt\big)^2}\\
                    &=& \frac{\big(\displaystyle\int_{I}g(\theta,t)dt\big)^2}{(n+1)^2 (n+1)^{\frac{n+1}{2}}\big(\displaystyle\int_{I}g(\theta,t)dt\big)^2}\\
                    &=& \frac{1}{(n+1)^2 (n+1)^{\frac{n+1}{2}}}.
\end{eqnarray*}
Since the lower bound obtained for $\alpha(n,\theta,I,K)$ is independent of $\theta$, $I$ and $K$, we have :
\begin{eqnarray*}
\alpha(n)\geq \frac{1}{(n+1)^2 (n+1)^{\frac{n+1}{2}}}.
\end{eqnarray*}

This proves Lemma \ref{nontriv}.

\end{proof}
\subsection{Main Theorem}
After that lengthy introduction, we are now ready to introduce the main theorem of this paper:
\begin{theorem} \label{main}
Let $n\geq 4$. Let $K$ be a symmetric convex set in $\mathbb{R}^n$. Then
\begin{eqnarray*}
vol_n(K)vol_n(K^{\circ})&\geq& \alpha(n-1)vol_{n-1}(\mathbb{S}^{n-1})^2\\
                    &=&\frac{4\alpha(n-1)\pi^n}{\Gamma(\frac{n}{2})^2},
\end{eqnarray*}
where $\alpha(n)$ is the constant defined in \ref{deff}.
\end{theorem} 

\begin{remark}
Applying Lemma \ref{nontriv} and Theorem \ref{main}, we obtain the following lower bound for the Mahler volume of every (at least four dimensional) symmetric convex set $K$ in $\mathbb{R}^n$ :
\begin{eqnarray*}
vol_n(K)vol_n(K^{\circ})\geq \frac{4\pi^n}{n^{\frac{n+4}{2}}\Gamma(\frac{n}{2})^2}.
\end{eqnarray*}
Obviously this lower bound is very uninteresting due to the fact that the lower bound in Lemma \ref{nontriv} is a very rough estimate.
\end{remark}
From this theorem, two different paths open up:
\begin{itemize}
\item First we attempt to calculate $\alpha(n)$ for $n\geq 4$ (for which I am less enthusiastic). 
\item Second, we attempt to prove that there exists a symmetric convex set $n\geq 4$ for which its Mahler volume is \emph{exactly} equal to the lower bound given by Theorem \ref{main} (which is more believable):
\end{itemize}

\begin{conjecture}\label{mine}
Let $n\geq 4$. There exists a symmetric convex set $K$ in $\mathbb{R}^n$ for which we have 
\begin{eqnarray*}
vol_n(K)vol_n(K^{\circ})=\frac{4\alpha(n-1)\pi^n}{\Gamma(\frac{n}{2})^2}.
\end{eqnarray*}
\end{conjecture}

\begin{remark}
\begin{itemize}
\item When the cone $C(I)$ is reduced to a point: the value of $\alpha(n,\theta,I,S)$ for every $S$ and every $\theta$ is equal to $1$ (i.e. set $\frac{0}{0}=1$)  and therefore:
\begin{eqnarray*}
\frac{4\alpha(n,\theta,I,S)\pi^n}{\Gamma(\frac{n}{2})^2}&=&\frac{4\pi^n}{\Gamma(\frac{n}{2})^2} \\
                                     &>& \frac{4^n}{\Gamma(n+1)}.
\end{eqnarray*}

\item When $C(I)$ is the half-plane: we shall discuss this case in the last section where an interesting lower bound for $\alpha(n,\theta,\mathbb{S}^1_{+},S)$ is conjectured.
\end{itemize}
\end{remark}

The rest of this paper is devoted to the proof of Theorem \ref{main}. To achieve this goal, we first require some techniques: in section $2$ we recollect the theory of convexly-derived measures, measurable partitions and the localisation on the canonical Riemannian sphere. Section $3$ deals with the proof of Theorem \ref{main}. And finally, Section $4$ contains further remarks concerning Theorem \ref{main} and Conjecture \ref{mine}.

\section{Localisation on the Sphere} \label{sin}

In the past few years, localisation methods have been used to prove several very interesting geometric inequalities. In \cite{lova} and \cite{kann}, the authors proved integral formulae using localisation, and applied their methods to conclude a few isoperimetric-type inequalities concerning convex sets in the Euclidean space. In \cite{guedon} the authors study a functional analysis version of the localisation, used again on the Euclidean space. Localisation on more general spaces was studied in \cite{gromil}, \cite{grwst}, \cite{memwst}, and \cite{memusphere}.

Many materials in this section are derived from \cite{memwst}. First, recall that a set $S\subset \mathbb{S}^n$ is called \emph{convex} (or geodesically convex) if for any $x,y\in S$, there exists a geodesic segment joining $x$ to $y$ which is contained in (the closure) of $S$.

\begin{definition}[Convexly-derived measures] \label{cderiv}
A convexly-derived measure on $\mathbb{S}^n$ (resp. $\mathbb{R}^n$) is a limit of a vaguely converging sequence of probability measures of the form $\mu_i=\frac{vol|S_i}{vol(S_i)}$, where $S_i$ are open convex sets. The space $\mathcal{MC}^n$ is defined to be the set of probability measures on $\mathbb{S}^n$ which are of the form $\mu_{S}=\frac{vol_{\vert S}}{vol(S)}$ where $S\subset \mathbb{S}^n$ is open and convex. The space of \emph{convexly-derived probability measures} on $\mathbb{S}^n$ is the closure of $\mathcal{MC}^n$ with respect to the vague (or weak by compacity of $\mathbb{S}^n$)-topology. The space $\mathcal{MC}^k$ will be the space of convexly-derived probability measures whose support has dimension $k$ and $\mathcal{MC}^{\leq k}=\cup_{l=0}^{k}\mathcal{MC}^{l}$. 
\end{definition}

This class of measures was defined first in \cite{gromil} and used later on in \cite{ale}, \cite{memwst}, \cite{memphd}. In Euclidean spaces, a convexly-derived measure is simply a probability measure supported on a convex set which has a $x^k$-concave density function with respect to the Lebesgue measure, for $k\geq 0$. A real function $f$ is $x^k$-concave if $f^{1/k}$ is a concave function. To understand convexly-derived measures on the sphere we need some definitions:

\begin{definition}[$\sin^k$-affine functions and measures]
A function $f$ is affinely $\sin^k$-concave if $f(x)=A\sin^k(x+x_0)$ for a $A>0$ and $0\leq x_0\leq \pi/2$. A $\sin^k$-affine measure by definition is a measure with a $\sin^k$-affine density function.
\end{definition}

\begin{lemma} \label{sinconc}
A real non-negative function $f$ defined on an interval of length less than $\pi$ is $\sin^k$-concave if for every $0<\alpha<1$ and for all $x_1,x_2 \in I$ we have
\begin{eqnarray*}
f^{1/k}(\alpha x_1+(1-\alpha)x_2)\geq (\frac{\sin(\alpha\vert x_2-x_1\vert)}{\sin(\vert x_2-x_1\vert)})f(x_1)^{1/k}+(\frac{\sin((1-\alpha)\vert x_2-x_1\vert)}{\sin(\vert x_2-x_1\vert)})f(x_2)^{1/k}.
\end{eqnarray*}
Particularly if $\alpha=\frac{1}{2}$ we have
\begin{eqnarray*}
f^{1/k}(\frac{x_1+x_2}{2})\geq \frac{f^{1/k}(x_1)+f^{1/k}(x_2)}{2\cos(\frac{\vert x_2-x_1\vert}{2})}.
\end{eqnarray*}
\end{lemma}
%\item $f$ admits only one maximum point and does not have any local minima.
%\item If $f$ is $\sin$-concave and defined on an interval containing $0$, then $g(t)=f(\vert t\vert)$ is also $\sin$-concave.
%\item Let $0<\varepsilon<\pi/2$. Let $\tau >\varepsilon$. $f$ is defined on $[0,\tau]$ and attains its maximum at $0$. Let $h(t)=c\cos(t)^k$ where $c$ is choosen such that $f(\varepsilon)=h(\varepsilon)$. Then
%\begin{eqnarray*}
%\begin{cases}
%f(x) \geq h(x)     & \text{for } x\in [0,\varepsilon], \\
%f(x) \leq h(x)     & \text{for } x\in [\varepsilon,\tau].
%\end{cases}
%\end{eqnarray*}
%In particular, $\tau\leq \pi/2$.
%\item Let $\tau>0$ and $f$ be a nonzero non-negative $\sin^k$-concave function on $[0,\tau]$ which attains its maximum at $0$. Then $\tau\leq \pi/2$ and for all $\alpha\geq 0$ and $\varepsilon\leq \pi/2$ we have
%\begin{eqnarray*}
%\frac{\displaystyle\int_{0}^{min\{\varepsilon,\tau\}}f(t)dt}{\displaystyle\int_{0}^{\tau}f(t)dt}\geq \frac{\displaystyle\int_{0}^{\varepsilon}\cos(t)^{k}dt}{\displaystyle\int_{0}^{\pi/2}\cos(t)^{k}dt}.
%\end{eqnarray*}
%\end{itemize}
%\end{lem}

This class of measures are also used in Optimal Transport Theory (see the excellent book \cite{villani}).

\begin{lemma} \label{memlem}
Let $S$ be a $k$-dimensional convex subset of the sphere $\mathbb{S}^n$ with $k \leq n$. Let $\mu$ be a convexly-derived measure defined on $S$ (with respect to the normalized Riemannian measure on the sphere). Then $\mu$ is a probability measure having a continuous density $f$ with respect to the canonical Riemannian measure on $\mathbb{S}^k$ restricted to $S$. Furthermore, the function $f$ is $\sin^{n-k}$-concave on every geodesic arc contained in $S$.
\end{lemma}
The above Lemma, proved in \cite{memwst}, completely characterizes the class of convexly-derived measures on the sphere. Note the similarity between the Euclidean case and the spherical one.

\begin{definition}[Spherical Needles] \label{sn}
A spherical needle in $\mathbb{S}^n$ is a couple $(I,\nu)$ where $I$ is a minimizing geodesic segment in $\mathbb{S}^n$ and $\nu$ is a probability measure \emph{supported} on $I$ which has a $\sin^{n-1}$-affine density function.
\end{definition}

\begin{remark}
\begin{itemize}
\item According to Lemma \ref{memlem}, the space of spherical needles is contained in $\mathcal{MC}^{\leq 1}$.
\item One can properly writes down the measure $\nu$. To do so, chose a parametrization of the geodesic segment $I$ by its arc length. Therefore there is a (canonical) map $s:[0,l(I)]\to I$. For every $t\in[0,l(s)]$, we have $\Vert \frac{ds}{dt}\Vert=1$. The measure $dt$ is the canonical Riemannian length-measure associated to the geodesic segment $I$. Then, $I$ is parametrized by $t\in [0,l(I)]$, the measure $\nu$ can be written as $\nu=C\cos(t-t_0)^{n-1}dt$, for $t_0\in[0,\pi]$ and $C$ is the normalization constant such that :
\begin{eqnarray*}
\int_{0}^{l(I)}C\cos(t-t_0)^{n-1}dt=1.
\end{eqnarray*}
\end{itemize}
\end{remark}
It is necessary to say a few words on \emph{convex partitions}. The reason being the fact that later on, one needs the canonical sphere (seen as a metric-measure space) to be partitioned into spherical needles. These objects need to be properly defined.
  
\begin{definition} \label{convexparti}
Let $\Pi$ be a finite convex partition of $\mathbb{S}^n$. We review this partition as an atomic probability measure $m(\Pi)$ on the space $\mathcal{MC}$ as follows: for each piece $S$ of $\Pi$, let $\mu_{S}=\frac{vol_{\vert S}}{vol(S)}$ be the normalised volume of $S$. Then set 
\begin{eqnarray*}
m(\Pi)=\sum_{S}\frac{vol(S)}{vol(\mathbb{S}^n)}\delta_{\mu_{S}}.
\end{eqnarray*}
Define the \emph{space} of (infinite) \emph{convex partitions} $\mathcal{CP}$ as the vague closure of the image of the map $m$ in the space $\mathcal{P}(\mathcal{MC})$ of probability measures on the space of convexly-derived measures. The subset $\mathcal{CP}^{\leq k}$ of convex partitions of dimension $\leq k$ consists of elements of $\mathcal{CP}$ which are supported on the subset $\mathcal{MC}^{\leq k}$ (of convexly derived measures with support of dimension at most $k$). It is worth remembering that the space $\mathcal{CP}$ is compact and $\mathcal{CP}^{\leq k}$ is closed within. 
\end{definition}

%We now include few useful deﬁnitions and remarks which we shall need  later on. The point being, the element of the partition of the sphere will always consist of convex sets, and there exist an algorithmic procedure which enables us to construct such sets. The following deﬁnitions are related to this fact

\begin{remark}
There exists an algorithmic procedure which enables one to construct the elements of the partition. The following definitions are related to this fact.
\end{remark}

\begin{definition}[Pancakes]
Let $S$ be an open convex subset of $\mathbb{S}^n$. Let $\varepsilon>0$. We call $S$ an $(k,\varepsilon)$-pancake if there exists a convex set $S_{\pi}$ of dimension $k$ such that every point of $S$ is at distance at most $\varepsilon$ from $S_{\pi}$. A set $S$ is called a \emph{pancake} if there exists a $k\geq 0$, $\varepsilon>0$ for which $S$ is a $(k,\varepsilon)$-pancake.
\end{definition}

\begin{remark} 
Every geodesic segment $I$ is a Hausdorff limit of a sequence $\{S_i\}_{i=1}^{\infty}$ , where $S_i$ is a $(1,\varepsilon_i)$-pancake and where $\varepsilon_i\to 0$ when $i\to\infty$. Furthermore, every spherical needle $(I,\nu)$ is a limit of a sequence of $(1,\varepsilon_i)$-pancakes where the measure $\nu$ is a (weak)-limit of the sequence of probability measures obtained by normalizing the volume of each pancake.
\end{remark}
\begin{definition}[Constructing Pancake]
Let $(I,\nu)$ be a spherical needle. We call a $(1,\delta)$-pancake $S$, a constructing pancake for $(I,\nu)$, if there exists a decreasing sequence of pancakes $...\subset S_i\subset S_{i-1}\subset\cdots\subset S_0$, where $S_0=S$ and $(I,\nu)$ is a limit of this sequence.
\end{definition}

The following Lemma will become useful in the next section, we will not present the proof and invite the reader to consult \cite{memwst}:
\begin{lemma} \label{restriction}
Let $\{S_i\}$ be a sequence of $(1,\delta_i)$-pancakes, where $\delta_i\to 0$ and the Hausdorff limit of the sequence $S_i$ is a geodesic segment $I$. The sequence of normalized measues $\mu_i=vol_{\vert S_i}/vol(S_i)$ vaguely converges to a convexly derived measure $\nu$ defined on $I$. Let $J=[x_1,x_2]\subset I$, where $[x_1,x_2]$ is the geodesic segment from the point $x_1$ to $x_2$. For $x\in J$, let $x^{\perp}$ denotes the $(n-1)$-dimensional sphere, $\mathbb{S}^{n-1}$, which contains $x$ and is \emph{orthogonal} to $J$. Then the sequence $\{A_i=S_i\cap x_1^{\perp}\cap x_2^{\perp}\}$ converges (in the Hausdorff metric) to the geodesic segment $J$ and the sequence of normalized measures $vol_{\vert A_i}/vol(A_i)$ vaguely converges to the meausre $\nu_{\vert J}/\nu(J)$.
\end{lemma}

\begin{definition}[Distance Between Spherical Needles]
For $\varepsilon>0$, we say that the distance between the spherical needles $(I_1,\nu_1)$ and $(I_2,\nu_2)$ is at most equal to $\varepsilon$ if there exists a constructing pancake $S_1$ (resp $S_2$) for $(I_1,\nu_1)$ (resp $(I_2,\nu_2)$) such that the Hausdorff distance (denoted by $d_H$) $d_H(S_1,S_2)$ is at most equal to $\varepsilon$.
\end{definition}
\begin{definition}[Cutting Map]
A \emph{cutting map} with respect to $(G_1,G_2,S)$, is a continuous map $f:\mathbb{S}^n\to\mathbb{R}^2$ defined by
\begin{equation} \label{eqn:cuttmap}
f(x)=\left(\int_{x^{\vee}\cap S}G_{1}(u)d\mu(u), \int_{x^{\vee}\cap S}G_{2}(u)d\mu(u)\right)
\end{equation}
where $x^{\vee}$ denotes the (oriented) open hemi-sphere centered at the point $x$, $G_i:\mathbb{S}^n\to\mathbb{R}$ (for $i=1,2$) are two continuous functions defined on the sphere, and $S$ is a \emph{convex} subset of $\mathbb{S}^n$. The normalized measure on the sphere is denoted by $d\mu$.
\end{definition}

\subsection{A Fundamental Spherical Localisation Lemma}
The main result of this section is the next lemma which is known as the \emph{localisation lemma}. The Euclidean counterpart is proved in \cite{lova}. The reader can skip the proof as it is very similar to the proof presented in \cite{lova}. 

\begin{lemma} \label{genlova}
Let $G_{i}$ for $i=1,2$  be two continuous functions on $\mathbb{S}^{n}$ such that
\begin{eqnarray*}
\int_{\mathbb{S}^{n}}G_{i}(u)d\mu(u)>0,
\end{eqnarray*}
then a convex partition of $\mathbb{S}^n$, $\Pi\in \mathcal{CP}^{\leq 1}$ by spherical needles exists such that for every $\sigma$ an element of $\Pi$, we have
\begin{equation} \label{eqn:baba}
\int_{\sigma}G_{i}(t)d\nu_{\sigma}(t)>0.
\end{equation}
$\nu_{\sigma}$ is a $\sin^{n-1}$-affine probability measure. 
\end{lemma}
%\begin{lem} \label{fondaa}
%Let $f_1$, $f_2$, $f_3$, and $f_4$ be four continuous non-negative functions on $\mathbb{S}^n$. Suppose that for every $u\in \mathbb{S}^n$ we have $f_1(u).f_2(u)\leq f_3(u).f_4(u)$. Let $\mu$ denote the normalised Riemannian measure of $\mathbb{S}^{n}$. If for every probability measure $\nu$ having a $\sin^{n-1}$-affine density function and supported on a geodesic segment $I$, we have
%\begin{eqnarray*}
%\left(\int_{I} f_1(t) d\nu(t)\right)\left(\int_{I} f_2(t) d\nu(t)\right)\leq \left(\int_{I} f_3(t) d\nu(t)\right)\left(\int_{I} f_4(t) d\nu(t)\right),
%\end{eqnarray*}
%where the geodesic segment $I$ is parametrised by its arc length $t$ and
%\begin{eqnarray*}
%f(t)d\nu(t)=f(t)g(t)dt,
%\end{eqnarray*}
%$g(t)$ is the density of the measure $\nu$, then we will have:
%\begin{equation} \label{eqn:gen}
%\left(\int_{\mathbb{S}^{n}} f_1(u) d\mu(u)\right)\left(\int_{\mathbb{S}^{n}}f_2(u) d\mu(u)\right)\leq \left(\int_{\mathbb{S}^{n}} f_3(u) d\mu(u)%%\right)\left(\int_{\mathbb{S}^{n}} f_4(u) d\mu(u)\right).
%\end{equation}
%\end{lem}

%This Lemma is a spherical version of the Localisation Lemma and its corollary which are both proved in \cite{lova} in the Euclidean case. I shall first prove the following:
%\begin{lem} \label{lova}
%Let $G_{i}$ for $i=1,2$  be two continuous functions on $\mathbb{S}^{n}$ such that
%\begin{eqnarray*}
%\int_{\mathbb{S}^{n}}G_{i}(u)d\mu(u)>0,
%\end{eqnarray*}
%then a $\sin^{n-1}$-affine probability measure $\nu$ supported on a geodesic segment $I$ exists such that
%\begin{eqnarray*}
%\int_{I}G_{i}(t)d\nu(t)>0.
%\end{eqnarray*}
%\end{lem}

\begin{proof}

First step is to prove the following claim:

\emph{claim}:

There exists a spherical needle $(I,\nu)$ such that equation (\ref{eqn:baba}) is satisfied.

\begin{proof}

We construct a decreasing sequence of convex subsets of $\mathbb{S}^n$ using the following procedure:
\begin{itemize}
\item Define the first step cutting map $F_1:\mathbb{S}^n \to \mathbb{R}^2$ by
\begin{equation} \label{eqn:borsuk}
F_1(x)=\left(\int_{x^{\vee}}G_{1}(u)d\mu(u), \int_{x^{\vee}}G_{2}(u)d\mu(u)\right)
\end{equation}
with respect to $(G_1,G_2,\mathbb{S}^n)$. Since \emph{cutting maps} are continuous, one can apply the Borsuk-Ulam Theorem to $F_1$. Hence there exists a $x_1^{\vee}$ such that
\begin{eqnarray*}
\int_{x_1^{\vee}}G_{1}(u)d\mu(u)=\int_{-x_1^{\vee}}G_{1}(u)d\mu(u) \\
\int_{x_1^{\vee}}G_{2}(u)d\mu(u)=\int_{-x_1^{\vee}}G_{2}(u)d\mu(u).
\end{eqnarray*}

Chose the hemi-sphere, denoted by $x^{\vee}_1$. Set $S_1=x^{\vee}_1\cap \mathbb{S}^{n}$.
\item Define the $i$-th step cutting map by
\begin{eqnarray*}
F_i(x)=\left(\int_{S_{i-1}\cap x^{\vee}}G_{1}(u)d\mu(u), \int_{S_{i-1}\cap x^{\vee}}G_{2}(u)d\mu(u)\right),
\end{eqnarray*}
 with respect to $(G_1,G_2,S_{i-1})$, where $S_{i-1}$ is the \emph{convex set} chosen in the $(i-1)$-th step. By applying the Borsuk-Ulam Theorem to $F_i$, we obtain two new hemi-spheres and we choose the one, denoted by $x^{\vee}_i$. Set $S_i=x^{\vee}_i\cap S_{i-1}$.
\end{itemize}

This procedure defines a decreasing sequence of convex subsets $S_i=x^{\vee}_i\cap S_{i-1}$ for every $i\in\mathbb{N}$. Set: 
 
\begin{eqnarray*}
S_{\pi}=\bigcap_{i=1}^{\infty}(S_{i})=\bigcap_{i=1}^{\infty}clos(S_{i}),
\end{eqnarray*}
where $clos(A)$ determines the topological closure of the subset $A$.
\begin{definition}[Cutting Hemi-spheres]
A \emph{solution of the Borsuk-Ulam Theorem} for a continuous map $f:\mathbb{S}^n\to\mathbb{R}^2$ is a point $x_0\in \mathbb{S}^n$ such that $f(x)=f(-x)$.

A cutting hemi-sphere is a $\mathbb{S}^n_{+}$ (half-sphere), which is a \emph{solution} of applying the Borsuk-Ulam Theorem to a \emph{cutting map}. The first cutting hemi-sphere will be the hemi-sphere used at the very first stage of the procedure to cut the sphere $\mathbb{S}^n$ into two parts.
\end{definition}

A convexly-derived probability measure $\nu_{\pi}$ is defined on $S_{\pi}$.
Since  $\lim_{i\to \infty}S_i=S_{\pi}$ (this limit is with respect to Hausdorff topology) the definition of the convexly-derived measures can be applied to define the positive probability measure supported on $S_{\pi}$ by
\begin{eqnarray*}
\nu_{\pi}=\lim_{i\to\infty}\frac{\mu|S_i}{\mu(S_i)}.
\end{eqnarray*}
Hence, by the definition of $\nu_{\pi}$
\begin{eqnarray*}
\int_{S_{\pi}}G_{j}(x)d\nu_{\pi}(x)=\lim_{i\to\infty}\frac{\displaystyle\int_{S_{i}}G_{j}(x)d\mu(x)}{\mu(S_{i})}
\end{eqnarray*}
for $j=1,2$, and where the limit is taken with respect to the vague topology defined on the space of convexly-derived measures (see \cite{memwst}). Recall the following :

\begin{lemma}
\label{radon}
{\em (See \cite{hirsh}).}
Let $\mu_i$ be a sequence of positive Radon measures on a locally-compact space $X$ which vaguely converges to a positive Radon measure $\mu$. Then, for every relatively compact subset $A \subset X$, such that $\mu(\partial A)=0$,
$$\lim_{i \to \infty}\mu_i (A)=\mu(A).$$
\end{lemma}

By the definition of the cutting maps $F_i(x)$, for every $i\in \mathbb{N}$, $j=1,2$ we have
\begin{eqnarray*}
\int_{S_{i}}G_{j}(u)d\mu(u)>0.
\end{eqnarray*}
By applying Lemma \ref{radon}, we conclude that the convexly-derived probability measure defined on $S_{\pi}$ satisfies the conclusion of the Lemma \ref{genlova}. The dimension of $S_{\pi}$ is $<n$. Indeed, if it is not the case, then $dim S_{\pi}=n$. Since there is a convexly-derived measure with positive density defined on $S_{\pi}$, and by the construction of the sequence $\{S_i\}$ for every open set $U$, there exists a $c>0$ such that we have
\begin{eqnarray*}
\nu_{\pi}(S_{\pi}\cap U)&=&\lim_{i\to\infty}\frac{\mu(S_i\cap U)}{\mu(S_i)}\\
                          &\leq&\lim_{i\to\infty}\frac{\mu(S_{\pi}\cap U)}{c^i\mu(S_{\pi})}.
\end{eqnarray*}
By supposition on the dimension of $S_{\pi}$, the right-hand side equality is equal to zero. This is a contradiction with the positive measure $\nu_{\pi}$ charging mass on $S_{\pi}\cap U$.

Thus, $dim S_{\pi} <n$. If $dim S_{\pi}\leq1$ then the claim is proved. If $dim S_{\pi}=k>2$, we define a new procedure by replacing $S$ with $S_{\pi}\cap S$, replacing the normalized Riemannian measure by the measure $\nu_{\pi}$, and replacing the sphere $\mathbb{S}^n$ by the sphere $\mathbb{S}^k$ containing $S_{\pi}$. For this new procedure, we define new cutting maps in every step. Since $k> 2$, by using the Borsuk-Ulam Theorem we obtain $\mathbb{S}^{k-1}$ halving the desired (convexly-derived) measures. The new procedure defines a new sequence of convex subsets and, by the same arguments given before, a convexly-derived measure defined on the intersection of this new sequence satisfying the conclusion of Lemma \ref{genlova}. By the same argument, the dimension of the intersection of the decreasing sequence of convex sets is $<k$. If the dimension of the intersection of this new sequence is at most equal to $1$,  we are finished. If not, we repeat the above procedure until arriving to an at most $1$-dimensional set. This proves that a probability measure $\nu$ with a (non-negative) $\sin^{n-1}$-concave density function $f$, supported on a geodesic segment $I$ exists such that:
\begin{equation} \label{eqn:dead}
\int_{I}f(t)G_i(t)dt\geq 0.
\end{equation}
We determine $I$ to have minimal length. If $f$ is $\sin^{n-1}$-affine on $I$ then we are done. We suppose this is not the case. We choose a sub-interval $J\subset I$, maximal in length, such that a $\sin^{n-1}$-concave function $f$ satisfying $(\ref{eqn:dead})$ exists such that $f$ additionally is $\sin^{n-1}$-affine on the sub-interval $J$. The existence of $J$ and $f$ follows from a standard compactness argument. We can assume that the length of $I$ is $<\pi/2$. Consider the Euclidean cone over $I$. Let $a,b\in I$ be the end points of $I$ and take the \emph{Euclidean segment} $[a,b]$ in $\mathbb{R}^2$ (basically the straight line joining $a$ to $b$). By definition of $\sin^{n-1}$-concave functions, the function $f$ is the restriction of a one-homogeneous $x^{n-1}$-concave function $F$ on the circle (a $x^{n-1}$-concave function $F$ is a function such that $F^{1/(n-1)}$ is concave). Transporting the entire problem to $\mathbb{R}^{n+1}$, we begin with two homogeneous functions $\bar{G_i}$ on $\mathbb{R}^{n+1}$ such that
\begin{eqnarray*}
\int_{\mathbb{R}^{n+1}}\bar{G_i} dx >0
\end{eqnarray*}
and we proved that there is a $2$-dimensional cone over a segment $[a,b]$, a one-homogeneous $x^{n-1}$-concave function $F$ on $[a,b]$, and a sub-interval $[\alpha,\beta]\subset [a,b]$ such that $F^{1/(n-1)}$ is linear on $[\alpha,\beta]$ (this is due to the fact that by definition, the restriction of a one-homogeneous $x^{n-1}$-affine function on a $2$-dimensional Euclidean cone defines a $\sin^{n-1}$-affine function on the circle) and such that
\begin{eqnarray*}
\int_{[a,b]}\bar{G_i}(t)F(t)dt\geq 0.
\end{eqnarray*}
We echo the arguments given in \cite{lova} (pages $21-23$) (with the only difference being that every construction there drops by one dimension). This drop of dimension is necessary so that that every construction may preserve homogeneity- or in other words, one dimension must be preserved for the $2$-dimensional cone defined on $[a,b]$. And the proof of the claim follows.
\end{proof}

In the proof of the claim we used a family of $\{x^{\vee}_i\}$ of oriented hemi-spheres to cut the sphere. Each $x^{\vee}_{i}$ cuts the sphere in two parts in such a way that in \emph{both} parts of the sphere the integral of $G_i$ remains positive (due to the Borsuk-Ulam Theorem). At each stage of the cutting, we only kept one part of the sphere. However, if we carry out everything we did with respect to the \emph{other} parts, we obtain (in a straightforward way) the conclusion of Lemma \ref{genlova}.

\end{proof}

\begin{remark}
\begin{itemize}
\item The proof of Lemma \ref{genlova} suggests that for any partition $\Pi\in\mathcal{CP}^{\leq 1}$ satisfying the conclusion of this lemma, there exists a sequence $\{\delta_i\}$, where for every $i\in\mathbb{N}$, $\delta_i>0$ and such that $\lim_{i\to\infty}\delta_i=0$, a sequence $\Pi_i\in\mathcal{CP}^{\leq n}$ such that the vague limit (in the sense of vague limit of measures) of the sequence $\Pi_i$ is equal to $\Pi$. Furthermore, for every $i\in\mathbb{N}$, every element of $\Pi_i$ is a $(1,\delta_i)$-pancake for (at least) a $1$-needle in $\Pi$.
\item We should keep in mind that the partition $\Pi\in \mathcal{CP}^{\leq 1}$ can be constructed by choosing the family of cutting hemi-spheres $\{x^{\vee}_i\}$ such that all the vectors $x_i$ belong to a sphere of dimension $k$ provided $k\geq 2$. This has two benefits:
\begin{itemize}
\item The partition obtained in Lemma \ref{genlova} is \emph{not} unique. Indeed we can define the \emph{cutting maps} to be defined on any $\mathbb{S}^k\subset \mathbb{S}^n$, provided $k\geq 2$.
\item We can choose the \emph{direction} of the cuts by appropriately choosing the sphere $\mathbb{S}^k$ from which we define the \emph{cutting maps}. This fact will become very useful in the proof of Theorem \ref{main}.
\end{itemize} 
\item Instead of applying the Borsuk-Ulam Theorem in the proof of Lemma \ref{genlova} we can use the more powerful Gromov-Borsuk-Ulam Theorem which is stated and proved in \cite{memwst}. The Gromov-Borsuk-Ulam Theorem provides a convex partition $\mathcal{CP}^{\leq k}$ for every continuous map $f:\mathbb{S}^n\to\mathbb{R}^k$ where $k<n$ and a point $z\in\mathbb{R}^k$ such that $f^{-1}(z)$ intersects the maximum points of the density of the convexly-derived-measures associated to the partition. Therefore by applying the Gromov-Borsuk-Ulam Theorem directly for the map $f:\mathbb{S}^n\to \mathbb{R}^2$ defined by
\begin{eqnarray*}
f(x)=\left(\int_{x^{\vee}}G_{1}(u)d\mu(u), \int_{x^{\vee}}G_{2}(u)d\mu(u)\right),
\end{eqnarray*}
we obtain the desired convex partition $\Pi\in\mathcal{CP}^{\leq 1}$ of the Lemma \ref{genlova}.
\end{itemize}
\end{remark}
It is interesting to compare the partition result of Lemma \ref{genlova} with the Yao-Yao partition (see \cite{yao} and \cite{leh2} for this matter).

\section{Proof of Theorem \ref{main}}

We are now ready to prove Theorem \ref{main}
\begin{proof}
The Mahler volume is invariant under linear mappings. According to John's Lemma (see \cite{john} or \cite{versh}), for every $n$-dimensional symmetric convex set $M$, there exists an ellipsoid of maximal volume $E$ contained in $M$, and such that $M$ is contained in $\sqrt{n}E$. Therefore, in order to calculate the Mahler volume of $M$, it is enough to calculate the Mahler volume of $L(M)$, where $L$ is an appropriate linear map and such that the ellipsoid of maximal volume contained in $L(M)$ coincides with $B_n(0,1)$. Hence, we will not lose any information by only evaluating the Mahler volume of the class of symmetric convex bodies in $\mathbb{R}^n$ for which their John's Ellipsoid is the unit ball $B_n(0,1)$. With this remark, we can then suppose $M$ is a symmetric convex body in $\mathbb{R}^n$ and
\begin{eqnarray*}
B_n(0,1)\subseteq M\subseteq B_n(0,\sqrt{n}).
\end{eqnarray*}

To prove Theorem \ref{main}, we shall proceed by contradiction. For simplicity we denote:
\begin{eqnarray*}
\beta(n)=\alpha(n-1)vol_{n-1}(\mathbb{S}^{n-1})^2.
\end{eqnarray*}
Assume we have a symmtric convex set in $\mathbb{R}^n$ such that:
\begin{equation} \label{eqn:onee}
vol_n(K)vol_n(K^{\circ})<\beta(n).
\end{equation}

Consider the normalised Riemannian volume on $\mathbb{S}^{n-1}$ which shall be denoted by $d\mu$. Therefore we have:
\begin{equation} \label{eqn:twoo}
\frac{vol_n(K)}{vol_{n-1}(\mathbb{S}^{n-1})}=\int_{\mathbb{S}^{n-1}}F_1(u)d\mu(u),
\end{equation}
where $F_1(u)=\frac{x(u)^n}{n}$ and $x(u)$ is the length of the segment issuing from the origin in the direction of $u$ and touching the boundary of $K$. Similarly, we have:
\begin{equation} \label{eqn:threee}
\frac{vol_n(K^{\circ})}{vol_{n-1}(\mathbb{S}^{n-1})}=\int_{\mathbb{S}^{n-1}}F_2(u)d\mu(u),
\end{equation}
where the function $F_2$ is defined similar to the function $F_1$.

%By assumption, we have:
%\begin{eqnarray*}
%\big(vol_{n-1}(\mathbb{S}^{n-1})\big)^2\big( \int_{\mathbb{S}^{n-1}}F_1(u)d\mu(u)\big)\big(\int_{\mathbb{S}^{n-1}}F_2(u)d\mu(u)\big)&<& \beta(n)\\
%                                                                                             &=&\alpha(n-1)\big(vol_{n-1}(\mathbb{S}^{n-1}\big)^2.
%\end{eqnarray*}
Combining equations (\ref{eqn:onee}), (\ref{eqn:twoo}) and (\ref{eqn:threee}), we obtain :
\begin{equation} \label{eqn:yes}
\big( \int_{\mathbb{S}^{n-1}}F_1(u)d\mu(u)\big)\big(\int_{\mathbb{S}^{n-1}}F_2(u)d\mu(u)\big)<\alpha(n-1).
\end{equation}

Thus, there exists a $C>0$ such that:
\begin{eqnarray*}
\int_{\mathbb{S}^{n-1}}F_1(u)d\mu(u) &<& C \\
                                     &<& \frac{\alpha(n-1)}{\displaystyle\int_{\mathbb{S}^{n-1}}F_2(u)d\mu(u)}.
\end{eqnarray*}

Define two functions $G_1$ and $G_2$ on $\mathbb{S}^{n-1}$ as follows. For every $u\in\mathbb{S}^{n-1}$, 
\begin{eqnarray*}
G_1(u)=C-F_1(u),
\end{eqnarray*}
and
\begin{eqnarray*}
G_2(u)=\alpha(n-1)-C\,F_2(u).
\end{eqnarray*}
Therefore, we have defined two continuous functions $G_1$ and $G_2$ on $\mathbb{S}^{n-1}$ such that (for every $i=1,2$) we have:
\begin{eqnarray*}
\int_{\mathbb{S}^{n-1}}G_i(u)d\mu(u)>0.
\end{eqnarray*}

We apply Lemma \ref{genlova} which demonstrates the existence of a partition of $\mathbb{S}^{n-1}$ into spherical needles, such that for every spherical needle $(I,\nu)$ we have:
\begin{eqnarray*}
\int_{I}G_i(t)d\nu(t)>0.
\end{eqnarray*}
This translates to the fact that for every such spherical needle $(I,\nu)$, we have :
\begin{equation} \label{eqn:darling}
\big(\int_{I}F_1(t)d\nu(t)\big)\big(\int_{I}F_2(t)d\nu(t)\big)< \alpha(n-1).
\end{equation}

We can observe here that by taking a $2$-dimensional (convex) cone over $I$ and writing the above integral inequality on this cone, we have translated the $n$-dimensional inequality (\ref{eqn:yes}) into a family of $2$-dimensional inequalities with respect to the measures $\mu_{2,\theta}$ of definition \ref{mes}.
\begin{remark}
From now on, any partition in $\mathcal{CP}^{\leq 1}$ is a partition $\Pi$ such that for any $1$-needle $(I,\nu)\in \Pi$, equation (\ref{eqn:darling}) is satisfied. The existence and construction of such partitions are provided by Lemma \ref{genlova}.
\end{remark}

The following definition could have been given in the introduction but we judged that introducing it here would be more meaningful:

\begin{definition}[Good Needles]
To every needle $(I,\nu)$ we assign the $\theta_0\in[0,\pi]$ which defines $\alpha((n-1),(L_I\cap M))$, where $L_I$ is the $2$-dimensional linear subspace of $\mathbb{R}^n$ which contains the convex cone $C(I)$.
Every needle $(J,\nu(\theta_0))$, where $J\subset(\mathbb{S}^{n-1}\cap L)$ is a geodesic segment of the unit circle, and $\nu(\theta_0)=g(t+\theta_0)dt$ is called a \emph{good needle}. The function $g$ is defined in definition \ref{deff}.
\end{definition}

\begin{lemma} \label{cut}
Let $n\geq 4$. For every $\varepsilon>0$, there exist a partition $\Pi$ and a spherical needle $(I_1,\nu_1)$ in $\Pi$ which is $\varepsilon$-close to a \emph{good needle} $(I_2,\nu_2)$.
\end{lemma}

\begin{proof}  

Let us define the following classes:
\begin{itemize}
\item The class $\mathcal{PA}$ is the class of partitions of $\mathbb{S}^{n-1}$ into $1$-needles such that every needle in this partition satisfies equation (\ref{eqn:darling}).
\item For every $\Pi\in\mathcal{PA}$, the class $\mathcal{NDL}(\Pi)$ is the class of $1$-needles belonging to the partition $\Pi$.
\item For every $1$-needle $(I,\nu)$, the class $\mathcal{PCK}(I,\nu)$ is the class of constructing pancakes for the needle $(I,\nu)$.
\item The class $\mathcal{GNDL}$ is the class of \emph{good needles}.
\item For every $(I,\nu)\in\mathcal{GPCK}$, the class $\mathcal{GPCK}(I,\nu)$ is the class of constructing pancakes for the needle $(I,\nu)$.
\end{itemize}

To prove Lemma \ref{cut}, we proceed by contradiction: 

%Let the $5$-tuples
%\begin{eqnarray*}
%(\Pi\in\mathcal{PA},(J,\nu_J)\in\mathcal{NDL}(\Pi),K_J\in\mathcal{PCK}(J,\nu_J),(G,\nu_G)\in\mathcal{GNDL},K_G\in\mathcal{GPCK}(G,\nu_G)).
%\end{eqnarray*}
We choose $(\Pi_0, (I_1,\nu_1), K_1,(I_2,\nu_2),K_2)$ such that
\begin{eqnarray*}
\Pi_0\in\mathcal{PA}\\
(I_1,\nu_1)\in\mathcal{NDL}(\Pi_0)\\
K_1\in\mathcal{PCK}(I_1,\nu_1)\\
(I_2,\nu_2)\in\mathcal{GNDL}\\
K_2\in\mathcal{GPCK}(I_2,\nu_2).
\end{eqnarray*}
 such that $d_H(K_1,K_2)$ is minimal over every $(\Pi,(I,\nu_I),K_I,(G,\nu_G),K_G)$ chosen as above. Let 
\begin{eqnarray*}
d_H(K_1,K_2)=\varepsilon_0>0.
\end{eqnarray*}

Since $\Pi_0$ is a partition of $\mathbb{S}^{n-1}$ into $1$-needles, therefore (according to the proof of Lemma \ref{genlova}) there exist a sequence $\delta_i\to 0$, a sequence $\Pi_i\in\mathcal{CP}^{\leq n}$ of partitions of $\mathbb{S}^{n-1}$ into $(1,\delta_i)$-pancakes such that the partition $\Pi_0$ is a (vague)-limit of the sequence $\Pi_i$. Therefore, there exists $i\in\mathbb{N}$ such that $K_1\in \Pi_i$ (this is exactly the definition of a \emph{constructing} pancake). From the partition $\Pi_i$, we like to obtain a \emph{new} limit partition, (perhaps) different from $\Pi_0$ (such that again, every needle in this new partition satisfies equation (\ref{eqn:darling})). In order to obtain this new partition, we shall cut each pancake $K_j\in\Pi_i$ applying a cutting map with respect to $(G_1,G_2,K_j)$. We concentrate only on cutting the \emph{constructing} pancake $K_1$ and we like the cuts to follow a particular direction: each cut should cut the pancake $K_1$ in such a way that the cutting hemi-sphere be orthogonal to the geodesic segment $I_2$. There exists a $(n-2)$- dimensional sphere such that each point of this sphere provides a cutting hemi-sphere orthogonal to $I_2$. Since $n\geq 4$, we can apply the Borsuk-Ulam theorem with respect to this $\mathbb{S}^{n-2}$. Therefore at each stage of the (algorithmic cutting) procedure, the cutting map will cut the pancake $K_1$ \emph{orthogonal} with respect to $I_2$. After the first cut, the pancake $K_1$ is divided in two parts (according to Borsuk-Ulam Theorem) which we denote by $K^1_1$ and $K^2_1$. The same scenario happens to the pancake $K_2$ \emph{i.e.} the pancake $K_2$ is also cut in two and we denote the two parts by $K^1_2$ and $K^2_2$. We choose the set $K^{i_0}_1$ where $i_0\in\{1,2\}$ such that for every $j\in\{1,2\}$ we have:
\begin{eqnarray*}
d_H(K^{i_0}_1,K^j_2)\leq d_H(K^s_1,K^p_2),
\end{eqnarray*}
for $s,p\in\{1,2\}$.
We denote the chosen set by $K^{i_0}_1$, which is the intersection of $K_1$ with the (first) cutting hemi-sphere. Then we start again cutting $K^{i_0}_1$ with the same strategy as before \emph{i.e.} cutting in a direction which is orthogonal to $I_2$ applying the Borsuk-Ulam theorem. After $N\geq 2$ \emph{finite} step cutting with the above strategy, we end up with a pancake $K^N_1\subset K_1$ and a pancake $K^N_2\subset K_2$ such that
\begin{equation} \label{eqn:ggg}
d_H(K^N_1,K^N_2)<d_H(K_1,K_2).
\end{equation}

Indeed this is the case since the \emph{longest} geodesic $\kappa$ from the boundary of $K_1$ orthogonal to the boundary of $K_2$ will eventually either be cut by the cutting hemi-spheres or the longest geodesic from the boundary of $K^N_1$ orthogonal to the boundary of $K^N_2$ will not contain $\kappa$, which concludes equation (\ref{eqn:ggg}).

The set $K^N_1$ is a constructing pancake for any needle $(W,\nu_W)$ in a limit partition, with $W\subset K^N_1$.

According to Lemma \ref{restriction}, the set $K^N_2$ is a constructing pancake for the good needle $(J,\nu_J)$, where $J=I_2\cap K^N_2$ and $\nu_J$ is the probability measure obtaining by restricting $\nu_2$ on $J\subset I_2$. Indeed this is the case since $K_2$ is consecutively cut in a direction orthogonal to $I_2$. 

According to equation (\ref{eqn:ggg}) we have :
\begin{eqnarray*}
d_H(K^N_1,K^N_2) &<& d_H(K_1,K_2)\\
                     &=& \varepsilon_0.
\end{eqnarray*}

This is a contradiction with the definition of $\varepsilon_0$ and the proof follows.

\end{proof}

Define a function 
\begin{eqnarray*}
F:\mathcal{MC}^{\leq 1}\to \mathbb{R}_{+},
\end{eqnarray*}
where $\mathcal{MC}^{\leq 1}$ is the space of \emph{spherical needles} by
\begin{eqnarray*}
F((I,\nu))=\big(\int_{I}F_1(t)d\nu(t)\big)\big(\int_{I}F_2(t)d\nu(t)\big).
\end{eqnarray*}

Since the space $\mathcal{MC}^{\leq 1}$ is compact, and since the function $F$ is a continuous function on $\mathcal{MC}^{\leq 1}$, for every $\varepsilon>0$, there exists a constant $C_1(\varepsilon)>0$ such that for every pair of spherical needles $(I_1,\nu_1)$ and $(I_2,\nu_2)$ which are $\varepsilon$-close to each other, we have :
\begin{equation} \label{eqn:continue}
\vert \big(\int_{I_1}F_1(t)d\nu_1(t)\big)\big(\int_{I_1}F_2(t)d\nu_1(t)\big)-\big(\int_{I_2}F_1(t)d\nu_2(t)\big)\big(\int_{I_2}F_2(t)d\nu_2(t)\big)\vert \leq C_1(\varepsilon).
\end{equation}

%Let $\{\varepsilon_i\}$ be a sequence of positive real numbers such that $\lim_{i\to\infty}\varepsilon_i=0$.

According to equations (\ref{eqn:darling})and (\ref{eqn:continue}), combined with Lemma \ref{cut}, there exists a good needle $(I_0,\nu_0)$ such that:
\begin{equation} \label{eqn:darling2}
\big(\int_{I_0}F_1(t)d\nu_0(t)\big)\big(\int_{I_0}F_2(t)d\nu_0(t)\big)< \alpha(n-1).
\end{equation}

%According to Lemma \ref{cut}, for every $\varepsilon>0$, there exists a $\Pi\in\mathcal{CP}^{\leq 1}$ and a spherical needle $(I,\nu)$ in $\Pi$ and a \emph{good needle} $(I_2,\nu_2)$ which are $\varepsilon$-close to each other.
Recall that for every $1$-needle $(I,\nu)$, we have:
\begin{equation} \label{eqn:tow1}
\big(\int_{I}F_1(t)d\nu(t)\big)\big(\int_{I}F_2(t)d\nu(t)\big)=\mu_{2,\theta}(C(I)\cap (M\cap L_I))\mu_{2,\theta}(C(I)\cap (M^{\circ}\cap L_I),
\end{equation}
where $L_I$ is the $2$-dimensional linear subspace of $\mathbb{R}^n$ which contains the convex cone $C(I)$ and $\mu_{2,\theta}$ is the measure as in definition \ref{deff} in which the function $g(t)$ coincides with the density of the probability measure $\nu$ of the $(I,\nu)$.

One part of the integral in equation (\ref{eqn:tow1}) deals with integrating over a $2$-dimensional section of $M^{\circ}$. By definition of the polar of a (symmetric) convex body we have:
\begin{eqnarray*} 
\int_{I}F_2(t)d\nu(t)&=&\mu_{2,\theta}(C(I)\cap (M^{\circ}\cap L_I)\\
                               &\geq&\mu_{2,\theta}(C(I)\cap (M\cap L_I)^{\circ}).
\end{eqnarray*}
According to the definition of the constant $\alpha(n)$ (definition \ref{deff}), for the \emph{good needle} $(I_0,\nu_0)$ we have:
\begin{eqnarray*}
\big(\int_{I_0}F_1(t)d\nu_0(t)\big)\big(\int_{I_0}F_2(t)d\nu_0(t)\big)\\
                                                              &=&\mu_{2,\theta_0}(C(I_0)\cap (L_{I_0}\cap M))\mu_{2,\theta_0}(C(I_0)\cap(L_{I_0}\cap M^{\circ}))\\
                                                              &\geq&\mu_{2,\theta_0}(C(I_0)\cap (L_{I_0}\cap M))\mu_{2,\theta_0}(C(I_0)\cap (L_{I_0}\cap M)^{\circ}) \\                                                              
                                                              &\geq& \alpha((n-1),(L_{I_0}\cap M))\\
                                                              &\geq&\alpha(n-1).
\end{eqnarray*}                                                              

Here $L_{I_{0}}$ is the $2$-dimensional linear subspace of $\mathbb{R}^n$ which contains the convex cone $C(I_0)$. This is a contradiction with equation (\ref{eqn:darling2}).                              
%To sum up, according to the above remarks and equation (\ref{eqn:tow1}), and lemma \ref{cut}, for every $\varepsilon>0$ we have constant $C'(\varepsilon)$ and a needle $(I_1,\nu_1)$ in a partition $\Pi$ and a \emph{good needle} $(I_2,\nu_2)$ such that:

%\begin{eqnarray*}
%\alpha(n-1)&>&\big(\int_{I_1}F_1(t)d\nu_1(t)\big)\big(\int_{I_1}F_2(t)d\nu_1(t)\big)\\
%&=&\mu_{2,\theta_1}(C(I_1)\cap (M\cap L))\mu_{2,\theta_1}(C(I_1)\cap (M^{\circ}\cap L)) \\
%&\geq& C'(\varepsilon)\big(\int_{I_2}F_1(t)d\nu_2(t)\big)\big(\int_{I_2}F_2(t)d\nu_2(t)\big)\\
%                                                              &=&\mu_{2,\theta_2}(C(I_2)\cap (L\cap M))\mu_{2,\theta_2}(C(I_2)\cap(L\cap M^{\circ}))\\
%                                                              &\geq&\mu_{2,\theta_2}(C(I_2)\cap (L\cap M))\mu_{2,\theta_2}(C(I_2)\cap (L\cap M)^{\circ}) \\                                                              
%                                                              &\geq& %C'(\varepsilon)\alpha(n-1).
%\end{eqnarray*}

%Which clearly is a contradiction as $\varepsilon\to 0$.

The proof of Theorem \ref{main} follows.

\end{proof}

\section{Remarks and Questions}

The first question which arises is how \emph{sharp} is Theorem \ref{main}? In order to demonstrate Conjecture \ref{mine}, one needs to study an inverse problem. Given the data $(n,\theta,I,S)$ in $\mathbb{R}^2$, can we construct an $n$-dimensional symmetric convex set $K$ such that:
\begin{equation} \label{eqn:no}
vol_n(K)vol_n(K^{\circ})=\alpha(n-1,\theta,I,S)vol_{n-1}(\mathbb{S}^{n-1})^2?
\end{equation}

%This seems plausible. Or perhaps, we could use a suitable probabilistic argument? Suppose the class of $n$-dimensional symmetric convex bodies is enhanced with a (suitable) probability measure. Perhaps we could show that the probability that a symmetric convex body (for which the equality (\ref{eqn:no}) holds) is equal to one? 

The second remark to point out is about a lower bound for $\alpha(n,\theta,\mathbb{S}^1_{+},S)$. The material for this follows from M.Fradelizi in \cite{fradel}. It seems reasonable to have: 
\begin{conjecture} \label{frad}
Let $\mathbb{R}^2$ be enhanced with the measure $\mu_2$ as in definition \ref{mes} and let $n\geq 4$. Let $K$ be a symmetric convex set in $\mathbb{R}^2$. Then we have:
\begin{eqnarray*}
\mu_2(K)\mu_2(K^{\circ})\geq \frac{1}{n}.
\end{eqnarray*}
\end{conjecture}

Conjecture \ref{frad} is valid for a certain class of symmetric convex bodies. Recall that an unconditional symmetric convex body in $\mathbb{R}^2$ is a convex body symmetric with respect to the origin of $\mathbb{R}^2$, which is also symmetric with respect to the $x$ and $y$-axis of the Cartesian coordinates. For unconditional symmetric convex bodies in $\mathbb{R}^2$, Conjecture \ref{frad} holds and can be obtained from the following result of Saint-Raymond:

\begin{theorem}[Saint-Raymond] \label{s}
Let $K$ be an unconditional symmetric convex body in $\mathbb{R}^n$. Let $K_{+}=K\cap\mathbb{R}^n_{+}$ and $K^{\circ}_{+}=K^{\circ}\cap\mathbb{R}^n_{+}$. Then for every $(m_1,\cdots,m_n)\in (0,\infty)^n$ we have:
\begin{eqnarray*}
\int_{K_{+}}\big(\prod_{i=1}^{n}x_i^{m_{i-1}}\big)dx\int_{K_{+}^{\circ}}\big(\prod_{i=1}^{n}x_i^{m_{i-1}}\big)dx\geq \frac{1}{\Gamma(m_1+\cdots+m_{n+1})}\prod_{i=1}^{n}\frac{\Gamma(m_i)}{m_i}.
\end{eqnarray*}
\end{theorem}

\begin{question}
Can one use symmetrization methods/shadow systems techniques similar to one used in \cite{fradel}, \cite{rog}, \cite{rog2} in order to prove Conjecture \ref{frad}?
This technique has shown interest in the study of Mahler Conjecture as one can consult \cite{fradshad}, \cite{fradpaou}, \cite{meyyr}.
Is it possible to generalize the proof of Theorem $1$ in \cite{camp} for our case (where the measure has an $\frac{1}{n}$-affine density $\vert x\vert^n$)?
\end{question}

Should Conjecture \ref{frad} be true, we will have:
\begin{eqnarray*}
\alpha(n,\theta,\mathbb{S}^1_{+},S)\geq \frac{1}{(n-1)C(n)^2},
\end{eqnarray*}
where 
\begin{eqnarray*}
C(n)=\int_{-\pi/2}^{+\pi/2}\cos(t)^{n-2}dt.
\end{eqnarray*}

And we will have:

\begin{eqnarray*}
\frac{vol_{n-1}(\mathbb{S}^{n-1})^2}{(n-1)C(n)^2}&=& \frac{vol_{n-2}(\mathbb{S}^{n-2})^2}{n-1} \\
                                                 &=& \frac{4\pi^{n-1}}{(n-1)\Gamma((n-1)/2)^2}\\
                                                 &>& \frac{4^n}{\Gamma(n+1)},
\end{eqnarray*}

which would suggest that the optimal value for $\alpha(n)$ is achieved on a convex cone which is strictly contained in the half-plane.

\bibliographystyle{plain}
\bibliography{mahl}

\begin{thebibliography}{10}

\bibitem{ale}
S.~{Alesker}.
\newblock Localization technique on the sphere and the {G}romov-{M}ilman
  theorem on the concentration phenomenon on uniformly convex sphere.
\newblock {\em Math. Sci. Res. Inst. Publ.}, 34:17--27, 1999.

\bibitem{shiri}
S.~Arstein-Avidan, R.~Karasev, and Y.~Ostrover.
\newblock From symplectic measurements to the {M}ahler {C}onjecture.
\newblock {\em arxiv:1303.4197}, 2013.

\bibitem{ball1}
K.~{Ball}.
\newblock {M}ahler's {C}onjecture and wavelets.
\newblock {\em Discrete Comput. Geom}, 13:271--277, 1995.

\bibitem{ball}
K.~Ball.
\newblock An elementary introduction to modern convex geometry.
\newblock In {\em Flavors of Geometry}, volume~31 of {\em MSRI lecture notes}.
  MSRI publications, 1997.

\bibitem{barthe}
F.~Barthe and M.~Fradelizi.
\newblock The volume product of convex bodies with many hyperplane symmetries.
\newblock {\em American Journal of Mathematics}, 135:1--37, 2013.

\bibitem{bl2}
W.~Blaschke.
\newblock Ueber affine geometrie: Neue extremeiegenschaften von ellipse und
  ellipsoid.
\newblock {\em Leipziger Ber.}, 69:306--318, 1917.

\bibitem{bl}
S.~G. {Bobkov} and M.~{Ledoux}.
\newblock Weighted {P}oincar\'e-type inequalities for {C}auchy and other convex
  measures.
\newblock {\em Ann. Probab.}, 37(2):403--427, 2009.

\bibitem{bourg}
J~Bourgain and V.D. Milman.
\newblock New volume ratio properties for convex symmetric bodies in
  $\mathbb{R}^n$.
\newblock {\em Invent. Math.}, 88:319--340, 1987.

\bibitem{camp}
S.~Campi and P.~Gronchi.
\newblock On volume product inequalities for convex sets.
\newblock {\em Proceedings of the American Mathematical society},
  134(8):2393--2402, 2006.

\bibitem{fradpaou}
D.~Cordero-Erausquin, M.~Fradelizi, P.~Paouris, and P.~Pivovarov.
\newblock volume of the polar of random sets and shadow systems.
\newblock {\em preprint}.

\bibitem{fradel}
M.~Fradelizi.
\newblock Private communication.
\newblock 2014.

\bibitem{guedon}
M.~{Fradelizi} and O.~{Gu\'edon}.
\newblock The extreme points of subsets of s-concave probabilities and a
  geometric localization theorem.
\newblock {\em Discrete Comput. Geom.}, 31:327--335, 2004.

\bibitem{fradmey}
M.~Fradelizi and M.~Meyer.
\newblock Functional inequalities related to {M}ahler's conjecture.
\newblock {\em Monatsh. Math.}, 159(1-2):13--25, 2010.

\bibitem{fradshad}
M.~Fradelizi, M.~Meyer, and A.~Zvavitch.
\newblock An application of shadow systems to {M}ahler's conjecture.
\newblock {\em Discrete and Computational Geometry}, 48(3):721--734, 2012.

\bibitem{gard}
R.~J. Gardner.
\newblock The {B}runn-{M}inkowski {I}nequality.
\newblock {\em Bull. Amer. Math. Soc.}, 39:355--405, 2002.

\bibitem{gord}
Y.~Gordon, M.~Meyer, and S.~Reisner.
\newblock Zonoids with minimal volume-product- a new proof.
\newblock {\em Proc. Amer. Math. Soc.}, 104(1):273--276, 1988.

\bibitem{groe}
H.~Groemer.
\newblock {\em Geometric applications of {F}ourier series and spherical
  harmonics}.
\newblock Cambridge University Press, 1996.

\bibitem{grwst}
M.~{Gromov}.
\newblock Isoperimetry of waists and concentration of maps.
\newblock {\em GAFA}, 13:178--215, 2003.

\bibitem{gromil}
M.~{Gromov} and V.D. {Milman}.
\newblock Generalisation of the spherical isoperimetric inequality to uniformly
  convex {B}anach spaces.
\newblock {\em Compositio Math.}, 62:3:263--282, 1987.

\bibitem{heuze}
M.~Henze.
\newblock The {M}ahler {C}onjecture.
\newblock {\em Thesis University of Magdeburg}, 2008.

\bibitem{hirsh}
F.~{Hirsch} and G.~{Lacombe}.
\newblock {\em Elements of Functional Analysis}, volume 192 of {\em Graduate
  Texts in Mathematics}.
\newblock Springer-Verlag, 1999.

\bibitem{hu}
P.~Hupp.
\newblock {M}ahler's {C}onjecture in convex geometry, a summary and further
  numerical analysis.
\newblock {\em Georgia Institute of Technology, Thesis}, 2010.

\bibitem{john}
F.~John.
\newblock Extremum problems with inequalities as subsidiary conditions.
\newblock {\em Interscience publishers Inc.}, 9:187--204, 1948.

\bibitem{kann}
R.~{Kannan}, L.~{Lov\'asz}, and M.~{Simonovits}.
\newblock Isoperimetric problems for convex bodies and a localization lemma.
\newblock {\em Discrete Comput. Geom.}, 13(3--4):541--559, 1995.

\bibitem{kuper1}
K.~Kuperberg.
\newblock A low-technology estimate in convex geometry.
\newblock {\em Internat. Math. Res. Notices.}, 9:181--183, 1992.

\bibitem{kuper}
K.~Kuperberg.
\newblock From the {M}ahler {C}onjecture to gauss linking integrals.
\newblock {\em Geom. Funct. Anal}, 18(3):870--892, 2008.

\bibitem{leh1}
J.~Lehec.
\newblock A direct proof of the functional {S}antalo inequality.
\newblock {\em C. R. Math. Acad. Sci. Paris}, 347(1-2):55--58, 2009.

\bibitem{leh2}
J.~Lehec.
\newblock On the {Y}ao-{Y}ao partition theorem.
\newblock {\em Arch. Math.}, 92(4):366--376, 2009.

\bibitem{lop}
M.~A. Lopez and S.~Reisner.
\newblock A special case of {M}ahler's {C}onjecture.
\newblock {\em Discrete Comput. Geom}, 20:163--177, 1998.

\bibitem{lova}
L.~{Lov\'asz} and M.~{Simonovits}.
\newblock Random walks in a convex body and an improved volume algorithm.
\newblock {\em Random {S}tructures {A}lgorithms}, 4(4):359--412, 1993.

\bibitem{lut}
E.~Lutwak and G.~Zhang.
\newblock Blaschke-santalo inequalities.
\newblock {\em J. Differential. Geom.}, 47:1--16, 1997.

\bibitem{reis2}
S.~M.~Reisner.
\newblock Zonoids with minimal volume product.
\newblock {\em Math. Z.}, 192:339--346, 1986.

\bibitem{mah1}
K.~Mahler.
\newblock Ein minimalproblem f\"ur konvexe polygone.
\newblock {\em Mathematica}, B7:118--127, 1938.

\bibitem{mah2}
K.~Mahler.
\newblock Ein \"ubertragungsprinzip f\"ur konvexe korper.
\newblock {\em Casopis Pest. Mat. Fys.}, 68:93--102, 1939.

\bibitem{memphd}
Y.~{Memarian}.
\newblock Geometry of the space of cycles: Waist and minimal graphs.
\newblock {\em Phd thesis}, 2010.

\bibitem{memwst}
Y.~{Memarian}.
\newblock On {G}romov's waist of the sphere theorem.
\newblock {\em Journal of {T}opology and {A}nalysis}, 3:7--36, 2011.

\bibitem{memusphere}
Y.~{Memarian}.
\newblock A lower bound on the waist of unit spheres of uniformly normed
  spaces.
\newblock {\em Compositio Math.}, 148(4):1238--1264, 2012.

\bibitem{meyer}
M.~Meyer.
\newblock Convex bodies with minimal volume product in $\mathbb{R}^2$.
\newblock {\em Monatsh. Math.}, 112:297--301, 1991.

\bibitem{meypaj}
M.~Meyer and A.~Pajor.
\newblock On the {B}laschke-{S}antal\'o inequality.
\newblock {\em Arch. Math.}, 55:82--93, 1990.

\bibitem{meyreis}
M.~Meyer and S.~Reisner.
\newblock Inequalities involving integrals of polar-conjugate concave
  functions.
\newblock {\em Monatsh. Math}, 125:219--227, 1998.

\bibitem{meyyr}
M~Meyer and S.~Reisner.
\newblock shadow systems and volume of polar convex bodies.
\newblock {\em Mathematika}, 53(1):129--148, 2007.

\bibitem{naz}
F.~Nazarov, F.~Petrov, D.~Ryabogin, and A.~Zvavitch.
\newblock A remark on the {M}ahler conjecture: local minimality of the unit
  cube.
\newblock {\em Duke. Math. J.}, 154(3):419--430, 2010.

\bibitem{pis}
G.~Pisier.
\newblock {\em The volume of convex bodies and {B}anach space geometry}.
\newblock Cambridge University Press, 1989.

\bibitem{reis}
S.~Reisner.
\newblock Random polytopes and the volume product of symmetric convex bodies.
\newblock {\em Math. Scand.}, 57:386--392, 1985.

\bibitem{reis3}
S.~Reisner.
\newblock Minimal volume-product in {B}anach spaces with a 1-unconditional
  basis.
\newblock {\em J. London. Math. Soc.}, 2(36):126--136, 1987.

\bibitem{reisut}
S.~Reisner, C.~Sch\"uett, and E.~M. Werner.
\newblock {M}ahler's conjecture and curvature.
\newblock {\em International Math. Research Notices}, 2012:1--16, 2012.

\bibitem{rog}
C.~A. Rogers and G.C. Shephard.
\newblock Some extremal problems for convex bodies.
\newblock {\em Mathematika}, 5:93--102, 1958.

\bibitem{sr}
J.~Saint-Raymond.
\newblock Sur le volume de corps convexes symetriques.
\newblock {\em Seminaire d'initiation a l'analyse}, 20, 1980/81.

\bibitem{santal}
L.~A. Santal\'o.
\newblock Un invariante afin para los cuerpos convexos del espacio de n
  dimensiones.
\newblock {\em Portugal. Math.}, 8:155--161, 1949.

\bibitem{sch}
R.~Schneider.
\newblock {\em Convex Bodies: The {B}runn-{M}inkowski theory}.
\newblock Cambridge University Press, 1993.

\bibitem{rog2}
G.~C. Shephard.
\newblock Shadow systems of convex bodies.
\newblock {\em Israel J. Math.}, 2:229--236, 1964.

\bibitem{stan}
A.~Stancu.
\newblock Two volume product inequalities and their applications.
\newblock {\em Canad. Math. Bull.}, 52:464--472, 2009.

\bibitem{tao2}
T.~Tao.
\newblock Open questions: the {M}ahler {C}onjecture on convex bodies.
\newblock {\em
  http://terrytao.wordpress.com/2007/03/08/open-problem-the-mahler-conjecture-on-convex-bodies/},
  2010.

\bibitem{tao1}
T.~Tao.
\newblock {S}antalo's inequality.
\newblock {\em www.math.ucla.edu/~tao/preprints/{E}xpository/santalo.dvi},
  2010.

\bibitem{versh}
R.~Vershynin.
\newblock Lectures in geometric functional analysis.
\newblock {\em Webpage www.personal.umich.edu/~romanv}.

\bibitem{villani}
C.~{Villani}.
\newblock {\em Optimal {T}ransport {O}ld {A}nd {N}ew}, volume 338.
\newblock Springer-Verlag, 2009.

\bibitem{white}
E.~White.
\newblock Polar-legendre duality in convex geometry and geometric flows.
\newblock {\em Georgia Institute of Technology. Master thesis in Mathematics}.

\bibitem{yao}
A.~C. Yao and F.~F. Yao.
\newblock A general approach to d-dimensional geometric queries.
\newblock {\em Proceedings of the seventeenth annual ACM symposium on Theory of
  computing,}, pages 163--168, 1985.

\end{thebibliography}

\end{document}